\documentclass[11pt]{amsart}

\usepackage{amsmath,amsthm,amsfonts,amssymb}
\usepackage{graphicx}
\usepackage{hyperref}
\usepackage[usenames]{color}

\theoremstyle{definition}

\newcounter{comcount}

\numberwithin{equation}{section}

\title{Cryptanalysis of Andrecut's public key cryptosystem}

\author{Vitali\u{i} Roman'kov}

\address{Institute of Mathematics and Information Technologies\\Omsk State Dostoevskii University}
\curraddr{}
\email{romankov48@mail.ru}

\author{Anton Menshov}

\address{Institute of Mathematics and Information Technologies\\Omsk State Dostoevskii University}
\curraddr{}
\email{menshov.a.v@gmail.com}
\thanks{}

\date{}

\begin{document}

\maketitle\footnote{Supported by RFBR, projects 13-01-00239 and 15-41-04312.}

\begin{abstract}
We show that a linear decomposition attack based on the decomposition method introduced by the first author in monography \cite{R1} and  papers \cite{R2}, \cite{MM}, and developed in \cite{RM}, works  by finding the exchanging key in the protocol in \cite{A}. 
\end{abstract}

\section{Introduction}
\label{se:intro}

In this note we apply  a practical  deterministic attack on the  protocol proposed in \cite{A}.
This kind of attack introduced by the first author in \cite{R1}, \cite{R2}, \cite{MM} and developed in \cite{RM} works when the platform  objects  are linear.
It turns out that in this case, contrary to the common opinion (and some explicitly stated security assumptions), one does not need to solve the underlying algorithmic problems to break the scheme, i.e., there is another  algorithm that recovers the private keys without solving the principal algorithmic problem on which the  security assumptions are based.
The efficacy of the attack depends on the platform group, so it requires a  specific analysis in each particular case.
In general one can only state that the attack is in polynomial time in the size of the data, when the platform and related groups are given together with their linear representations.
In many other cases we can effectively use known linear presentations of the groups under consideration.
A theoretical base for the decomposition method is described in \cite{RM} where a series of examples is presented.
The monography \cite{R1} solves uniformly  many protocols based on the conjugacy search problem, protocols based on the decomposition and factorization problems, protocols based on actions by automorphisms, and a number of other protocols.
See also \cite{R3} and \cite{R4} where the linear decomposition attack is applied to the  main protocols in \cite{WXLLW},  \cite{HKKS}, and \cite{KLS}.  

In a series of works \cite{T1}, \cite{T2} and \cite{T3} (see also \cite{BKT}) Tsaban presented another general approach for provable polynomial time solutions of computational problems in groups with efficient, faithful representation as matrix groups. 

All along the paper we denote by $\mathbb{N}$ the set of all nonnegative integers, and by $\mathbb{C}$ the set of all complex numbers.

\section{Andrecut's key exchange protocol \cite{A}.}

In this section, we describe the Andrecut's key exchange protocol proposed in \cite{A}.
Firstly we introduce a necessary terminology.
Then we will give a cryptanalysis of this protocol. 

Let $X \in \mathbb{C}^{n\times n}$ be a complex matrix of size $n \times n$, that is considered as a variable. Let 
\[
P(X,a) = \sum_{m=1}^M a_m X^m \  \textrm{and} \  
Q(X,b) = \sum_{m=1}^K b_m X^m
\]

\noindent
be two complex polynomials in $X$ uniquely defined by the complex vectors of coefficients $a = (a_0,\dots,a_{M})$ and $b= (b_0,\dots,b_{K}).$ 

\begin{itemize}
\item Alice chooses the secret vectors  $a \in \mathbb{C}^{M_1}$ and $\tilde{a} \in \mathbb{C}^{M_2}$ (Alice's private key).
\item Alice randomly generates and publishes the matrix $U\in \mathbb{C}^{n\times n}$ (Alice's matrix public key). 
\item Bob  chooses the secret vectors $b \in \mathbb{C}^{J_1}$ and $\tilde{b} \in \mathbb{C}^{J_1}$ (Bob's  private key).
\item Bob randomly generates and publishes the matrix $V\in \mathbb{C}^{n\times n}$ (Bob's matrix public key).
\item Alice computes and publishes the matrix $A = P(U,a)P(V,\tilde{a})$ (Alice's public key).
 \item Bob computes and publishes the matrix $B = P(U,b)P(V,\tilde{b})$ (Bob's public key).
\item Alice calculates the secret matrix $K_A = P(U,a) B P(V,\tilde{a}).$ 
\item Bob calculates the secret matrix $K_B = P(U,b)AP(V,\tilde{b})$. 
\item The established secret key is $K = K_A = K_B.$
\end{itemize}

It is assumed in \cite{A}, that the matrices $U$ and $V$ are different to give the non-commutativity assumption $P(U,a)P(V,b) \not= P(V,b)P(U,a)$ and  $P(U,\tilde{a})P(V,\tilde{b}) \not= P(V,\tilde{b})P(U,\tilde{a})$.
In general, even more strong assumption $UV\not= VU$ is not enough for this non-commutativity. 

Also, there is a remark in \cite{A} that the following assumption
\[
M_1, M_2, J_1, J_2 \gg n
\] 
should be satisfied in order to increase the security.
But by the classical Cayley-Hamilton theorem every matrix $U$ is a root of its own  characteristic polynomial $C(X) = \det(U - X\cdot I_n)$, where $I_n$ is the identity matrix of size $n \times n$.
The degree of $C(X)$ is exactly $n$.
Then for any matrix $U\in \mathbb{C}^{n\times n}$ and every matrix polynomial $P(X,a)$ one has $P(U,a)=R(U),$ where $R(X)$ is the remainder after division of the polynomial $P(X,a)$ by $C(X).$
Hence, there is no sense to use in the protocol above polynomials of degrees $\geq n.$ 

\section{Cryptanalysis of the Andrecut's key exchange protocol \cite{A}.}

We will provide two approaches to cryptanalysis of the protocol described in the previous section.
The first one is based on some simple facts from linear algebra and the second one is based on a linear decomposition attack.

As we have seen in the previous section, there is no sense to use values $M_1,M_2,J_1,J_2 \geq n$.
Thus we can assume that $M_1=M_2=J_1=J_2=n-1$.
Consider the linear systems
\begin{equation}\label{eq:systems}
\begin{aligned}
A &= \sum_{i=1}^{n-1} \sum_{j=1}^{n-1} \underbrace{a_i \tilde{a}_j}_{x_{ij}} U^i V^j, \\
B &= \sum_{k=1}^{n-1} \sum_{l=1}^{n-1} \underbrace{b_k \tilde{b}_l}_{y_{kl}} U^k V^l
\end{aligned}
\end{equation}
of $n^2$ equations with $(n-1)^2$ unknowns $x_{ij}$ and $y_{kl}$.
Having a solution of the systems (\ref{eq:systems}) one can compute the secret key as follows
\begin{equation}\label{eq:key_computation}
\begin{aligned}
K &= P(U,a) P(U,b) P(V,\tilde{a}) P(V,\tilde{b}) \\
  &= \left( \sum_{i=1}^{n-1} \sum_{j=1}^{n-1} a_i b_j U^{i+j} \right) \left( \sum_{k=1}^{n-1} \sum_{l=1}^{n-1} \tilde{a}_k \tilde{b}_l V^{k+l} \right) \\
  &= \sum_{i=1}^{n-1} \sum_{j=1}^{n-1} \sum_{k=1}^{n-1} \sum_{l=1}^{n-1} \underbrace{a_i\tilde{a}_k}_{x_{ik}} \underbrace{b_j\tilde{b}_l}_{y_{jl}} U^{i+j} V^{k+l}.
\end{aligned}
\end{equation}
Solution of a system of $n$ equations with $n$ unknowns using Gauss elimination requires $O(n^3)$ time.
Thus solution of the systems (\ref{eq:systems}) requires $O(n^6)$ time.
Having precomputed values $U^i$ and $V^i$, for $i=1,\dots,n-1$, one can perform the step (\ref{eq:key_computation}) in $O(n^7)$ time.
So the overall time complexity for this approach is $O(n^7)$.

Further we will describe the second approach based on a linear decomposition attack.
The algebra $\mathbb{C}^{n\times n}$ has a structure of a vector space over $\mathbb{C}$ of dimension $n^2.$
Let $\bar{U}$ be the semigroup generated by $U,$ and let 
$\bar{V}$ be the semigroup generated by $V.$
A basis of the subspace Sp$(\bar{U}\bar{V})$ can be effectively constructed as follows.
Let
\begin{align*}
L_0 &= \{ I_n \}, \\
L_1 &= \{ U, V \}, \\
    &\dots \\
L_i &= \{ U^k V^l \mid k,l \in \mathbb{N}, \ k+l = i \}, \\
    &\dots
\end{align*}
be the sets of matrices considered as vectors.
Define
\[
V_i=\mathrm{Sp}(L_0\cup L_1 \cup \dots \cup L_i),
\]
for $i=0,1,\dots,$ the vector space spanned by the indicated set.  
We choose a basis $B_0 = \{b_0 = I_n\}$ of $V_0$, then extend $B_0$ to basis $B_1$ of $V_1$, and so on.
If for some $i_0$ we get $B_{i_0}=B_{i_0+1},$ then clearly $B = B_{i_0}$ is a basis of Sp$(\bar{U}\bar{V})$.
By the Cayley-Hamilton theorem one has $i_0 \leq 2n-2.$
In construction of $B$  we only  use the Gauss elimination process that is polynomial.
Note, that we can do it offline, so we will call this phase the {\it offline phase}.
Let $b_0, b_1, \dots, b_r$ be a basis of Sp$(\bar{U}\bar{V})$, where $b_i = U^{k_i}V^{l_i}$.
Now we are ready to recover the secret key $K$ ({\it online phase}).

\begin{itemize}
\item Since $B\in $ Sp$(\bar{U}\bar{V})$ we can use the Gauss elimination process to obtain  a presentation of $B$ in the form
\begin{equation}\label{eq:B}
B = \sum_{i=0}^r \alpha_i U^{k_i} V^{l_i}, \quad \alpha_i \in \mathbb{C}.
\end{equation}
\item Then we have
\begin{equation}\label{eq:key_computation_2}
\begin{aligned}
\sum_{i=0}^r\alpha_iU^{k_i}AV^{l_i} &=
\sum_{i=0}^r\alpha_iU^{k_i}P(U,a)P(V,\tilde{a})V^{l_i} \\ &=
P(U,a) \left( \sum_{i=0}^r\alpha_iU^{k_i}V^{l_i} \right) P(V,\tilde{a}) \\ &=
P(U,a)BP(V,\tilde{a}) = K.
\end{aligned}
\end{equation}
\end{itemize}
  
Note that a similar protocol by Stickel \cite{S} has been analyzed in \cite{RM}.
A linear decomposition attack based on the decomposition method has been applied.

Now we will provide a rough estimate for the time complexity of the approach above.
Observe that the number of the field operations in Gauss elimination performed on a matrix of size $k \times n^2$ is $O(k^2n^2)$.
A basis of Sp$(\bar{U}\bar{V})$ consist of at most $n^2$ elements, so it requires $O(n^2 \sum_{k=1}^{n^2} k^2)=O(n^8)$ time to construct it.
Computing $\alpha_i$ in (\ref{eq:B}) by solving a system of linear equations using Gauss elimination requires $O(n^6)$ time.
Having precomputed values $U^{k_i}$ and $V^{l_i}$, for $i=0,\dots,r$, one can perform the step (\ref{eq:key_computation_2}) in $O(n^5)$ time.
So the offline phase could be done in $O(n^8)$ time, the online phase could be done in $O(n^6)$ time, and the overall time complexity is $O(n^8)$.

\end{document}